\documentclass[11pt,a4paper]{article}

\usepackage[english]{babel}
\usepackage[latin1]{inputenc}
\usepackage[T1]{fontenc}
\usepackage{graphicx}
\usepackage{amsmath}
\usepackage{amssymb}
\usepackage{latexsym}
\usepackage{amsthm}
\usepackage[all]{xy}
\usepackage{stmaryrd}
\usepackage{subfigure}
\usepackage{MnSymbol}

\def\be{\begin{equation}}
\def\ee{\end{equation}}
\def\beq{\begin{eqnarray*}}
\def\eeq{\end{eqnarray*}}
\def\rk{\operatorname{rk}}
\def\Ker{\operatorname{Ker}}
\def\Im{\operatorname{Im}}
\def\Z{\mathbb{Z}}

\newcommand{\pic}[3]{\parbox[c]{#1cm}{$\includegraphics[scale=#2]{#3}$}}

\newtheorem{teo}{Theorem}[section]
\newtheorem{lem}{Lemma}[section]

\newtheorem{prop}{Proposition}[section]
\theoremstyle{definition}
\newtheorem{defn}{Definition}[section]
\theoremstyle{remark}

\begin{document}

\begin{center}
{\large \textbf{Step by step categorification of the Jones polynomial in Kauffman's version} } \\
Alessio Carrega \\
Department of Mathematics ``L. Tonelli'',\\
University of Pisa,\\
Largo Bruno Pontecorvo 5 56127 Pisa, Italy\\
carrega@mail.dm.unipi.it
\end{center}

\subsection*{Abstract}
Given any  diagram of a link, we define on the cube of Kauffman's states a ``$2$-complex'' whose homology is an invariant of the associated framed links, and such that the graded Euler characteristic reproduces the unnormalized Kauffman bracket. This includes a categorification of brackets skein relation. Then we incorporate the orientation information and get a further complex on the same cube that gives rise to a new invariant homology for oriented links, so that the graded Euler characteristic reproduces the unnormalized Jones polynomial in Kauffman's version. Finally we clarify the relations between this homology and the original Khovanov homology of oriented links, extending the well known relation between the associated two versions of the Jones  polynomial.

\subsection*{Keywords}
Kauffman bracket, Jones polynomial, frame, Khovanov homology, framed links, computable, invariant, oriented links

\section{Introduction}

Kauffman's derivation \cite{Kauffman} of the Jones polynomial of an oriented link, starts with the brackets state sum over any diagram of the link (non oriented and equipped with the black-board framing) interpreted as an unnormalized invariant of the framed link; then modifies the brackets by forgetting the framing and incorporating the orientation information carried by the writhe, eventually obtaining the unnormalized Jones polynomial in Kauffman's version.

The Khovanov homology of an oriented link \cite{Khovanov} is obtained from a complex of graded $\Z$-modules supported by the cube of Kauffman states of any non oriented diagram of the link. The orientation information is incorporated applying some shift. The graded Euler characteristic of this homology reproduces the unnormalized Jones polynomial in Khovanov's version. There is a well known easy relation ($q=-A^{-2}$) between these two versions of the Jones polynomial of an oriented link.

The aim of this paper is to adapt Khovanov's constructions and ideas (mostly following Bar-Natan's algebraic exposition \cite{Bar-Natan}), in order to provide a step by step categorification of Kauffman's derivation. For any diagram of a link (non oriented and equipped with the black-board framing), with the same cube of Kauffman's states as support, we define a new ``$2$-complex'' whose homology is an invariant of the associated framed link. The graded Euler characteristic reproduces the unnormalized Kauffman bracket. This includes a categorification of the skein relations that determines the brackets. Then, by forgetting the framing and incorporating the orientation information, we produce a further complex whose homology is a invariant of the oriented link, so that the graded Euler characteristic reproduces the unnormalized Jones polynomial in Kauffman's version.

Although this homology of oriented links is formally new, the respective constructions are so close to each other, that one expects a strict relation with the original Khovanov homology, extending the relation between the two versions of the Jones polynomial. These relations are completely clarified at the end of the paper. We note that existing efficient algorithms in order to compute Khovanov homology should be easily adapted to compute the homologies defined in this paper (in particular the framed link one). We conclude this paper showing some examples and proving that our homology for framed links is actually a generalization of Kauffman bracket, namely there exist framed links (in particular we show framed knots) with the same Kauffman bracket but different homology.

Khovanov's categorification is based on Khovanovs's version of the (unnormalized) Jones polynomial, in particular on its writing as sum on the Kauffman states. The Khovanov's version of the Jones polynomial is useful for obtain a generalization because, from the second equation of the skein relations we have a power of the polynomial $q+q^{-1}$ on each summand of the Jones polynomial, and is easy to find a graded abelian group whose graded Euler characteristic is $q+q^{-1}$. Kauffman's version is not so good since we have the polynomial $-A^2 - A^{-2}$ instead of $q+q^{-1}$, and is impossible to find a graded abelian group whose graded Euler characteristic is $-A^2-A^{-2}$. One of the problems that we had in make our construction was precisely ``how to put the less inside the powers of $A^2+A^{-2}$''. As we said, at the end of this paper we show the connection between the classical Khovanov homology and our invariant for oriented links. This result says  that Khovanov homology generalize our construction, but this is not surprising since our idea for obtain a correct proof of the invariance was to merge  the even components and the odd ones of the complex.

\section{Preliminary algebraic definitions}

\begin{defn}
Given $i\in\{0,1\}$ we denote with $\underline{i}$ the element of $\{0,1\}$ different from $i$.
\end{defn}
\begin{defn}
Let $\mathcal{C}$ be a preadditive category. $2\!-\!\mathcal{C}_*(\mathcal{C})$ is the category defined by the followings:
\begin{itemize}
\item{\textit{objects}: pairs of objects of $\mathcal{C}$, $(X_0,X_1)$, together with two arrows of $\mathcal{C}$, $\partial_0 : X_0 \rightarrow X_1$ e $\partial_1 : X_1 \rightarrow X_0$, such taht $\partial_1 \circ \partial_1 =0$ and $\partial_0 \circ \partial_1= 0$. Those maps are called \emph{differentials} or \emph{boundary homomorphisms} respectively of degree $0$ e $1$. We'll write $X=(X_0 ,X_1 , \partial_0, \partial_1)$. These objects are said to be $2$-\emph{complexes} in $\mathcal{C}$.}
\item{\textit{morphisms}: pairs of morphisms of $\mathcal{C}$, $(f_0: X_0 \rightarrow Y_0 , f_1: X_1 \rightarrow Y_1)$, that commute with the differentials: $f_{\underline{i}} \circ \partial_i = \partial_i \circ f_i$ for each $i \in\{0,1 \}$.}
\item{\textit{identities}: pairs of identities of $\mathcal{C}$.}
\item{\textit{composizione}: component whise.}
\end{itemize}
\end{defn}
We are interested in $2\!-\!\mathcal{C}_*(\mathcal{G}\textit{r}\mathcal{A}\textit{b})$ and $2\!-\!\mathcal{C}_*(2\!-\!\mathcal{C}_*(\mathcal{G}\textit{r}\mathcal{A}\textit{b}))$, where $\mathcal{G}\textit{r}\mathcal{A}\textit{b}$ is the category of graded abelians groups and morphisms of graded abelian groups of degree $0$.

\begin{defn}
 A $2$-\emph{subcomplex} of a $2$-complex of graded abelian groups $A=(A_0, A_1, \partial_0, \partial_1)$ is a $2$-complex $A'=(A'_0 , A'_1, \partial'_0, \partial'_1)$ such that for each $i\in\{0,1\}$ $A'_i$ is a graded subgroup of $A_i$ and $\partial'_i$ is obtained by restriction of $\partial_i$.
\end{defn}
\begin{defn}
Let $A'=(A'_0 , A'_1, \partial'_0, \partial'_1)$ be a $2$-subcomplex of a $2$-complex of abelian groups $A=(A_0, A_1, \partial_0, \partial_1)$. The \emph{quotient} of $A$ modulo $A'$ is the $2$-complex defined by $A/A' :=(A_0/A'_0, A_1/A'_1, \bar \partial_0, \bar \partial_1)$ where the differentials are given by passing to the quotient of abelian groups.
\end{defn}
In the same way we can define the subobjects and the quotients in the case of an arbitrary abelian category, in particular for $2\!-\!\mathcal{C}_*(\mathcal{G}\textit{r}\mathcal{A}\textit{b})$.
\begin{defn}
Given two $2$-coplexes in an additive category $\mathcal{C}$, $A=(A_0, A_1, \partial^A_0, \partial^A_1)$ and $B=(B_0, B_1, \partial^B_0, \partial^B_1)$, we define the \emph{direct sum} of $A$ and $B$ as the $2$-complex $A\oplus B :=(A_0\oplus B_0, A_1\oplus B_1, \partial^A_0\oplus \partial^B_0 , \partial^A_1\oplus \partial^B_1)$.
\end{defn}
With these notions the category of $2$-complexes in an additive category becomes additive. If moreover $\mathcal{C}$ is an abelian category, we can define the kernels and the images of morphisms in $2\!-\!\mathcal{C}_*(\mathcal{C})$, obtaining an abelian category.

\begin{defn}
The \emph{homology} of a $2$-complex of graded abelian groups $A=(A_0, A_1, \partial_0, \partial_1)$ is the pair of graded abelian groups $(H_o(A), H_1(A))$ defined by:
$$
H_i(A) := \frac{\Ker \partial_i}{ \Im \partial_{\underline{i}} }
$$
Given a morphism of $2$-complexes $f=(f_0,f_1): A \rightarrow B$, with $A=(A_0, A_1, \partial^A_0, \partial^A_1)$ and $B=(B_0, B_1, \partial^B_0, \partial^B_1)$, this induces a map in homology for each $i\in\{0,i\}$ $f_{*,i} : H_i(A) \rightarrow H_i(B)$ obtained as a restriction to the kernels  and passing at the quotients.\\
Given $i\in\{0,1\}$ we define the following functor, called $i$-th \emph{homology functor}:
$$
H_i: 2\!-\!\mathcal{C}_*(\mathcal{G}\textit{r}\mathcal{A}\textit{b} ) \longrightarrow \mathcal{G}\textit{r}\mathcal{A}\textit{b}
$$
$$
\xymatrix{
A \ar[dd]_{f} \ar@{|->}[rr] & & H_i(A) \ar[dd]^{f_{*,i}} \\
 \ar@{|->}[rr] & & \\
B \ar@{|->}[rr] & & H_i(B)
}
$$
Likewise we can define the homology functor on the category of $2$-complexes in an arbitrary abelian category, in particular for $2\!-\!\mathcal{C}_*(\mathcal{G}\textit{r}\mathcal{A}\textit{b})$.
\end{defn}

\begin{defn}
Given a $2$-complex $X=(X_0, X_1, \partial_0, \partial_1)$ we define the \emph{reflexed} of $X$: $X^\spadesuit :=(X_1, X_0, \partial_1, \partial_0)$. The \emph{reflexed} of a map of $2$-complexes $f: X\rightarrow Y$ is the morphism $f^\spadesuit :=(f_1,f_0) : A^\spadesuit \rightarrow B^\spadesuit$.
\end{defn}

\begin{defn}
Given a $2$-complex of $2$-complex of graded abelian groups $A=(A_0, A_1, \partial_0, \partial_1)$, with $\partial_i = (\partial_i^0, \partial_i^1)$ and $A_i=(A_{i,0}, A_{i,1}, \partial_{i,0}, \partial_{i,1})$ for each $i\in\{0,1\}$, the \emph{flatten} of $A$ is the $2$-complex of graded abelian groups $\textit{Fl} (A) := ( A_{0,0} \oplus A_{1,1}, A_{0,1}\oplus A_{1,0}, (\partial_{0,0} \oplus \partial_{1,1}) + (\partial_0^1\oplus \partial_1^0) , (\partial_{0,1}\oplus \partial_{1,0}) + (\partial_0^0\oplus \partial_1^1))  $. Hence at the level of spaces $\textit{Fl}(A) = A_0 \oplus A^\spadesuit_1$.\\
Given a morphism of $2$-complexes of $2$-complexes of graded abelian groups $f=(f_0,f_1): A \rightarrow B$, with $A=(A_0, A_1, \partial^A_0, \partial^A_1)$, $B=(B_0, B_1, \partial^B_0, \partial^B_1)$, $A_i=(A_{i,0}, A_{i,1}, \partial^A_{i,0}, \partial^A_{i,1})$, $B_i=(B_{i,0}, B_{i,1}, \partial^A_{i,0}, \partial^B_{i,1})$ and $f_i = (f_{i,0} ,f_{i,1})$ for each $i\in\{0,1\}$, we define the morphism of $2$-complexes $\textit{Fl}(f) := (f_0 \oplus f^\spadesuit_1) : \textit{Fl}(A) \rightarrow \textit{Fl} (B)$.\\
We define the \emph{flatting functor}:
$$
\textit{Fl}: 2\!-\!\mathcal{C}_*(2\!-\!\mathcal{C}_* ( \mathcal{G}\textit{r}\mathcal{A}\textit{b} ) ) \longrightarrow 2\!-\!\mathcal{C}_* ( \mathcal{G}\textit{r}\mathcal{A}\textit{b} )
$$
$$
\xymatrix{
A \ar[dd]_{f} \ar@{|->}[rr] & & \textit{Fl}(A) \ar[dd]^{\textit{Fl}(f)} \\
 \ar@{|->}[rr] & & \\
B \ar@{|->}[rr] & & \textit{Fl}(B)
}
$$
\end{defn}
We have that for each $i\in\{0,1\}$
$$
\textit{Fl}\circ H_i = H_i \circ \textit{Fl} : 2\!-\!\mathcal{C}_*(2\!-\!\mathcal{C}_* ( \mathcal{G}\textit{r}\mathcal{A}\textit{b} ) ) \rightarrow \mathcal{G}\textit{r}\mathcal{A}\textit{b}
$$

\begin{defn}
Let $A=(A_0, A_1, \partial_0, \partial_1)$ be a $2$-complex of graded abelian groups. The \emph{graded Euler characteristic} of $A$ is the Laurent polynomial
$$
\chi_A (A) := q\!\dim A_0 - q\!\dim A_1
$$
\end{defn}

\begin{defn}
Let $A=(A_0, A_1, \partial^A_0, \partial^A_1)$ and $B=(B_0, B_1, \partial^B_0, \partial^B_1)$ be two $2$-complexes of graded abelian groups. The \emph{tensor product} of $A$ and $B$ is the $2$-complex $A\otimes B :=((A_0\otimes B_0) \oplus (A_1\otimes B_1), (A_0\otimes B_1) \oplus (A_1\otimes B_0), [\partial^A_0 \otimes \textit{id}_{B_0} + \textit{id}_{A_0} \partial^B_0, \partial^A_1 \otimes \textit{id}_{B_1} - \textit{id}_{A_1} \partial^B_1] , [\partial^A_0 \otimes \textit{id}_{B_1} + \textit{id}_{A_0} \partial^B_1, \partial^A_1 \otimes \textit{id}_{B_0} - \textit{id}_{A_1} \partial^B_0] )$.
\end{defn}
It holds that the graded Euler characteristic of the direct sum of two $2$-complexes is the sum of the characteristics, that is
$$
\chi_A (A \oplus B) = \chi_A (A) + \chi_A(B)
$$
The characteristic of the tensor product is the product of the characteristics, that is
$$
\chi_A (A \otimes B) = \chi_A (A) \cdot\chi_A(B)
$$

\begin{teo}
Let
$$
0 \longrightarrow A \stackrel{f}{\longrightarrow} B \stackrel{g}{\longrightarrow} C \longrightarrow 0
$$
be a short sequence in $2\!-\!\mathcal{C}_*( \mathcal{C} )$, with $\mathcal{C}$ an abelian category. Then we have an homology exact sequence
$$
\xymatrix{
H_0(A) \ar[r]^{f_{*,0}} & H_0(B) \ar[r]^{g_{*,0}} & H_0(C) \ar[d]^{\Delta_0} \\
H_1(C) \ar[u]^{\Delta_1} & H_1(B) \ar[l]^{g_{*,1}} & H_1(A) \ar[l]^{f_{*,1}}
}
$$
\begin{proof}
As in the classic case.
\end{proof}
\end{teo}
So we have a version of the lemma used by Bar-Natan for proving the invariance of Khovanov homology for oriented links in \cite{Bar-Natan}:

\begin{lem}\label{lemquoziente}
Let $A$ be a $2$-complex in an abelian category and $A'$ a subcomplex. Then
\begin{itemize}
\item{if $A'$ has homology $0$, then $H^*(A/ A') \cong H^*(A)$;}
\item{if $A/A'$ has homology $0$, then $H^*(A) \cong H^*(A')$.}
\end{itemize}
\begin{proof}
It follows from the long exact sequence induced by the short exact sequence
$$
\begin{matrix}
0 & \rightarrow & A' & \stackrel{j}{\hookrightarrow } & A & \stackrel{\pi}{\rightarrow} & \frac{A}{A'} & \rightarrow & 0
\end{matrix}
$$
\end{proof}
\end{lem}
We remark that this lemma has not to be used in the category $2\!-\!\mathcal{C}_*(\mathcal{G}\textit{r}\mathcal{A}\textit{b})$, but in $2\!-\!\mathcal{C}_* ( 2\!-\!\mathcal{C}_* ( \mathcal{G}\textit{r}\mathcal{A}\textit{b} ) )$.

\section{Costruction of the invariant}

\begin{defn}
We define the graded abelian group $W$ as the group with each homogeneous component $0$ except in degree $2$ end $-2$, where it has the free abelian groups $W_+$ and $W_-$ respectively generated by $w_+$ e $w_-$.
$$
W: \begin{matrix}
\ldots & 0 & W_- & 0 & 0 & 0 & W_+ & 0 & \ldots \\
\ldots & -3 & -2 & -1 & 0 & 1 & 2 & 3 & \ldots
\end{matrix}
$$
\end{defn}

\begin{defn}
Given a link diagram $D$ and a Kauffman state $s$ of $D$, we define
$$
W_s(D) := \left( \bigotimes_{\pic{0.2}{0.1}{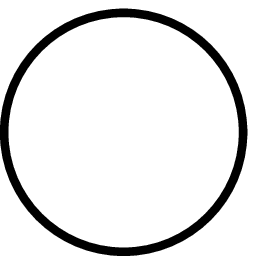} \text{ in } D_s} W \right) \{a(s) -b(s)\}
$$
Where $D_s$ is the splitting of $D$ by the state $s$, $a(s)$ is the number of crossings of $D$ marked with $A$ (or $0$) by $s$, and $b(s)$ is the number of crossings marked with $B$ (or $1$).
\end{defn}

\begin{defn}
Let $D$ be a link diagram. For each $i\in\{0,1\}$ we define the following graded abelian group
$$
\llangle D \rrangle_i := \bigoplus_{s\ : \ b(s) + |s_A| \in 2\Z + i} W_s(D)
$$
where $s_A$ is the Kauffman state of $D$ maid only by $A$ (or only $0$), and $|s|$, for a state $s$ of $D$, is the number of components of the splitting of $D$ by $s$.
\end{defn}
\begin{defn}
We define the following maps of graded abelian groups of degree $2$ similarly to the classic $m$ and $\Delta$ by replacing $v_+$ with $w_-$ and $v_-$ with $w_+$:
$$
\begin{matrix}
\begin{array}{rcl}
W\otimes W & \stackrel{\bar m}{\longrightarrow} & W \\
w_- \otimes w_- & \longmapsto & w_- \\
w_+ \otimes w_- & \longmapsto & w_+ \\
w_- \otimes w_+ & \longmapsto & w_+ \\
w_+ \otimes w_+ & \longmapsto & 0
\end{array}

&

\begin{array}{rcl}
W & \stackrel{\bar \Delta}{\longrightarrow} & W \otimes W \\
w_- & \longmapsto & w_+\otimes w_- + w_-\otimes w_+ \\
w_+ & \longmapsto & w_+\otimes w_+
\end{array}
\end{matrix}
$$
\end{defn}

\begin{defn}
Given an edge $\xi:s_0 \rightarrow s_1$ of a diagram $D$, we define the homomorphism of graded abelian groups of degree $0$
$$
\partial_\xi : W_{s_0}(D) \longrightarrow W_{s_1}(D)
$$
likewise the classic $d_\xi : V_{s_0}(D) \rightarrow V_{s_1}(D)$ using $\bar m$ and $\bar \Delta$ instead of $m$ and $\Delta$
$$
\partial_\xi := \left\{ \begin{array}{cl}
\{-2\} ( \textit{id}_W \otimes \ldots \otimes \textit{id}_W \otimes \bar m \textit{id}_W \otimes \ldots \otimes \textit{id}_W ) \{ a(s_0) - b(s_0) \} & \text{if} \ |s_1| = |s_0| + 1 \\
\{-2\} ( \textit{id}_W \otimes \ldots \otimes \textit{id}_W \otimes \bar \Delta \textit{id}_W \otimes \ldots \otimes \textit{id}_W ) \{ a(s_0) - b(s_0) \} & \text{if} \ |s_1| = |s_0| - 1
\end{array}\right.
$$
\end{defn}
The maps defined above are of degree $0$. in fact $\bar m$ and $\bar \Delta $ are of degree $2$, $a(s_1) = a(s_0) - 1 $ and $b(s_1)=b(s_0) +1$, hence $a(s_1) - b(s_1) = a(s_0) - b(s_0) - 2$. By definition of $W_s(D)$, we obtain that the degree of $\partial_\xi$ is $0$.

\begin{defn}
Let $D$ be a link diagram, together with an order of its crossings. Mimicking what we have done in the classic case to define the differentials $d^j : \mathcal{C}^j(D) \rightarrow \mathcal{C}^{j+1}(D)$ with $j\in\Z$, for each $i\in\{ 0, 1\}$ we define the map
$$
\partial_i := \sum_{\xi \ :\ |\xi| + |s_A| \in 2\Z + i} (-1)^\xi \partial_\xi : \llangle D \rrangle_i \longrightarrow \llangle D \rrangle_{\underline{i}}
$$
\end{defn}
The definition is well done. In fact, given an edge $\xi: s_0 \rightarrow s_1$, $\partial_\xi$ is of degree $0$, $b(s_1)=b(s_0)+1$, hence $b(s_1)+|s_A| \in 2\Z \Leftrightarrow b(s_0)+|s_A| \in 2\Z +1$.

\begin{lem}
Let $D$ be a link diagram with an order of the crossings and let $i \in\{0,1\}$. Then
$$
\partial_{\underline{i}} \circ \partial_i = 0
$$
\begin{proof}
It follows from the commutativity of the induced cube, such as in the classic case.
\end{proof}
\end{lem}

So we can define the following $2$-complexes of graded abelian groups:
\begin{defn}
Let $D$ be a link diagram with an enumeration of the crossings
$$
\llangle D \rrangle := ( \llangle D \rrangle_0, \llangle D \rrangle_1 , \partial_0 , \partial_1)
$$
\end{defn}

\begin{teo}\label{carEul}
Let $D$ be a link diagram. Then, using the variable $A$
$$
\chi_A( \llangle D \rrangle ) = (-A^{-2} - A^2) \langle D \rangle
$$
Thus the graded Euler characteristic of $\llangle D \rrangle $ is equal to the unnormalized Kauffman bracket of $D$.
\begin{proof}
We remember that $s_A$ is the state of $D$ with all $A$.
\beq
\chi_A(\llangle D \rrangle ) & = & q\!\dim \llangle D \rrangle_0 - q\!\dim \llangle D \rrangle_1 \\
 & = & \sum_{s\ : \ b(s) + |s_A| \in 2\Z} q\!\dim W_s(D) - \sum_{s\ : \ b(s) + |s_A| \in 2\Z + 1 } q\!\dim W_s(D) \\
 & = & \sum_{s \text{ of }D } (-1)^{b(s)+|s_A|} q\!\dim W_s(D) \\
 & = & \sum_{s \text{ of }D } (-1)^{b(s)+|s_A|} A^{a(s)- b(s)} (q\!\dim W )^{ |s|} \\
 & = & \sum_{s \text{ of }D } (-1)^{b(s)+|s_A|} A^{a(s)-b(s)} (A^{-2} + A^2)^{|s|}
\eeq
Given an edge $\xi:s_0 \rightarrow s_1$, we have that $b(s_1)=b(s_0)+1$ and $|s_1| = |s_0| \pm 1$. Each state $s$ can be linked to $s_A$ with a sequence of edges whose length $b(s)$
$$
s_A \stackrel{\xi_1}{\longrightarrow} \ldots \stackrel{\xi_{b(s)}}{\longrightarrow} s
$$
Hence the class modulo $2$ of $b(s)$ is the same as $|s|- |s_A|$, therefore $(-1)^{b(s)+|s_A|} = (-1)^{|s|}$. Hence
\beq
\chi_A(\llangle D \rrangle ) & = & \sum_{s \text{ of }D } (-1)^{|s|} A^{a(s)-b(s)} (A^{-2} + A^2)^{|s|} \\
 & = & \sum_{s \text{ of }D } A^{a(s)-b(s)} (-A^{-2} - A^2)^{|s|} \\
 & = & (-A^{-2}-A^2) \langle D \rangle
\eeq
\end{proof}
\end{teo}

\section{Proof of the invariance}
We denote with $\llangle \ \rrangle $ the application defined on the set of the link diagrams with an order of the crossings and values in the class of $2$-complexes of graded abelian groups that associates to a diagram $D$ the $2$-complex $\llangle D \rrangle$. We denote with $H_*: 2\!-\!\mathcal{C}_*( \mathcal{G}\textit{r}\mathcal{A}\textit{b} ) \rightarrow (\mathcal{G}\textit{r}\mathcal{A}\textit{b})^2$ the homology functor joining the functors $H_0$ and $H_1$.
\begin{teo}
$H_* \circ \llangle \ \rrangle $ is an invariant for framed links, up to isomophism.
\end{teo}
The isomorphism for two different orders of crossings can be obtained as for the Khovanov homology for oriented links (see \cite{Lee}). The invariance for the change of a curl with another of the same type is trivial. Following the proofs of the invariance for the Reidemeister moves of the second and third type for Khovanov homology in \cite{Bar-Natan} we obtain proofs valid also for this version. In this section we adapt the Bar-Natan's proof about the moves of the second type to this case. Doing the same adaptations for the moves of the third type we conclude.\\

We denote with $FL(\mathbb{S}^3)$ the set of the framed links of the $3$-sphere up to equivalence. Thanks to the theorem, we can give the following definition:
\begin{defn}
$$
\mathcal{H}^F_* : FL(\mathbb{S}^3) \longrightarrow \{ \text{Pairs of graded abelian groups} \}/_{\cong}
$$
$$
\mathcal{H}^F_*(L) := [H_*( \llangle D \rrangle ) ]/_{\cong}
$$
 where $D$ is a diagram of the framed link $L$ with a fixed enumeration of the crossings. Given $L\in FL(\mathbb{S}^3)$, $i\in\{0,1\}$ and $j\in\Z$ we denote with $\mathcal{H}^F_i(L)$ the $i$-th component of $\mathcal{H}^F_*(L)$, and with $\mathcal{H}^F_{i,j}(L)$ the $j$-th homogeneous component of $\mathcal{H}^F_i(L)$. Given a diagram $D$ with an order of the crossings, we use $\mathcal{H}^F_*(D)$ instead of $H_*(\llangle D \rrangle )$ and in analogous way for the components. If the diagram has not a fixed order of the crossings, we use the same notation to indicate the item up to isomorphism.
\end{defn}

\begin{teo}
Let $L$ be a framed link. Then the graded Euler characteristic of $\mathcal{H}^F_*(L)$ is equal to the Kauffman bracket of $L$.
$$
\chi_A(\mathcal{H}^F_*(L)) = \langle L \rangle
$$
\begin{proof}
It follows from the Theorem \ref{carEul} and from the invariance of the graded Euler characteristic by passing in homology.
\end{proof}
\end{teo}

As usual we implicitly give an order on the first presented diagram and we consider the order induced from that one to the others. Given a diagram with a local behavior of the type $\pic{0.9}{0.2}{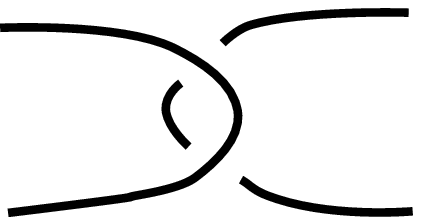}$, we have the following decomposition in direct sum for the spaces:
$$
\left\llangle \pic{1.2}{0.3}{reid2-1p.eps} \right\rrangle = \left\llangle \pic{1.2}{0.3}{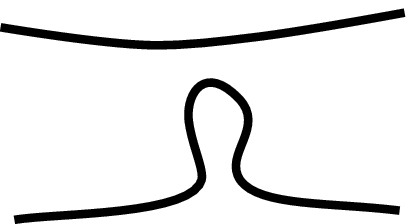} \right\rrangle\{2\} \oplus \left\llangle \pic{1.2}{0.3}{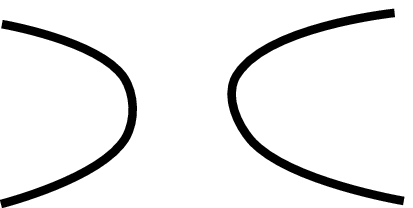} \right\rrangle \oplus \left\llangle \pic{1.2}{0.3}{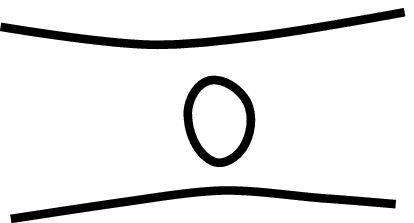} \right\rrangle \oplus \left\llangle \pic{1.2}{0.3}{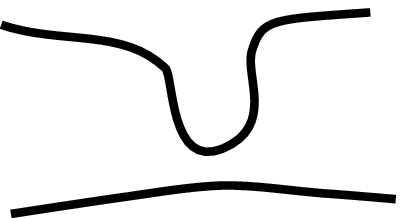} \right\rrangle \{-2\}
$$
We explain the decomposition. Let $s$ be a Kauffman state of one of the diagrams in the right member of the equation. We denote with $s_A$ the state of the diagram with all $A$ (or $0$), and we denote with $s'_A$ the state of \pic{0.9}{0.2}{reid2-1p.eps} with all $A$. $s$ corresponds to a state $s'$ of \pic{0.9}{0.2}{reid2-1p.eps}.
\begin{itemize}
\item{If $s$ is of \pic{0.9}{0.2}{reid2c.eps} then $a(s)=a(s')-2$, $b(s)=b(s')$, $a(s)-b(s)=a(s')-b(s')-2$, hence the polynomial degree has to be shifted by $2$. $|s_A| = |s'_A|$, hence $b(s) + |s_A| \in 2\Z \Leftrightarrow b(s') + |s'_A|$, it is correct not to reflex.}
\item{If $s$ is of \pic{0.9}{0.2}{Bcanalep.eps} then $a(s)=a(s')-1$, $b(s)=b(s')-1$, $a(s)-b(s)=a(s')-b(s')$, hence the polynomial degree has to not be shifted. $|s_A| = |s'_A| \pm 1$, depend if the pieces of the link that we can see in the figure \pic{0.9}{0.2}{reid2-1p.eps} stay in the same component or not, therefore $b(s) + |s_A| \in 2\Z \Leftrightarrow b(s') + |s'_A|$, it is correct not to reflex.}
\item{If $s$ is of \pic{0.9}{0.2}{reid2d.eps} then $a(s)=a(s')-1$, $b(s)=b(s')-1$, $a(s)-b(s)=a(s')-b(s')$, hence the polynomial degree has not to be shifted. $|s_A| = |s'_A| + 1$, hence $b(s) + |s_A| \in 2\Z \Leftrightarrow b(s') + |s'_A|$, it is correct not to reflex.}
\item{If $s$ is of \pic{0.9}{0.2}{reid2e.eps} then $a(s)=a(s')$, $b(s)=b(s')-2$, $a(s)-b(s)=a(s')-b(s')+2$, hence the polynomial degree has to be shifted by $-2$. $|s_A| = |s'_A|$, hence $b(s) + |s_A| \in 2\Z \Leftrightarrow b(s') + |s'_A|$, it is correct not to reflex.}
\end{itemize}
We define the $2$-complex of $2$-complexes $C$ such that $\textit{Fl}(C) = \left\llangle \pic{0.9}{0.2}{reid2-1p.eps} \right\rrangle$:
$$
\begin{matrix}
C:= \left( \left\llangle \pic{1.2}{0.3}{Bcanalep.eps} \right\rrangle \oplus \left\llangle \pic{1.2}{0.3}{reid2d.eps} \right\rrangle , \left\llangle \pic{1.2}{0.3}{reid2c.eps} \right\rrangle^\spadesuit\{2\} \oplus \left\llangle \pic{1.2}{0.3}{reid2e.eps} \right\rrangle^\spadesuit\{-2\} ,\right.\\
, [( 0, \partial_{B*} ), (0, \bar m )], [(\partial_{*A}, \partial_{A*} ) , 0  ]  \Bigr{ ) }
\end{matrix}
$$
where the maps that compose the differentials are built summing the maps induced from the edges of \pic{0.9}{0.2}{reid2-1p.eps} that change the crossings in figure multiplied by the sign of the edge. These maps are pieces of the differentials of $\left\llangle \pic{0.9}{0.2}{reid2-1p.eps} \right\rrangle$, therefore taken $i\in\{0,1\}$ we have that the $i$-th component of these goes from the component of degree $i$ to the one of degree $\underline{i}$, therefore if we reflex the codomain or the domain we obtain the respect of the degree.\\
We define the $2$-subcomplex of $C$, \\
$$
C':= \left(\left\llangle \pic{1.2}{0.3}{reid2d.eps} \right\rrangle_- , \left\llangle \pic{1.2}{0.3}{reid2e.eps} \right\rrangle^\spadesuit \{-2\}, \bar m , 0 \right)
$$
where $\left\llangle \pic{0.9}{0.2}{reid2d.eps} \right\rrangle_-$ is the $2$-subcomplex of $\left\llangle \pic{0.9}{0.2}{reid2d.eps} \right\rrangle$ obtained considering only the component $W_-$ of the factor of the tensor products that correspond to the circle in figure. The differentiale of degree $0$ is given by $\bar m$ restricted to the $2$-subcomplex. We can easily see that it is an isomorphism, and therefore that $C'$ has homology $0$. Hence the homology of $C$ is isomorphic to the homology of
$$
\frac{C}{C'} = \left( \left\llangle \pic{1.2}{0.3}{Bcanalep.eps} \right\rrangle \oplus \left\llangle \pic{1.2}{0.3}{reid2d.eps} \right\rrangle_+ , \left\llangle \pic{1.2}{0.3}{reid2c.eps} \right\rrangle^\spadesuit\{2\} , 0 , [(\partial_{*A}, \partial_{A*} ) , 0  ] \right)
$$
The map $\partial_{A*}: \left\llangle \pic{0.9}{0.2}{reid2c.eps} \right\rrangle \{2\} \rightarrow \left\llangle \pic{0.9}{0.2}{reid2d.eps} \right\rrangle $ is an isomorphism. We can define the map $\tau:= \partial_{*A} \circ \partial_{A*}^{-1} : \left\llangle \pic{0.9}{0.2}{reid2d.eps} \right\rrangle \rightarrow \left\llangle \pic{1.2}{0.3}{Bcanalep.eps} \right\rrangle $ and the $2$-subcomplex of $C/C'$,
$$
C''' := \left( \left\{ (\tau(b), b) \in \left\llangle \pic{1.2}{0.3}{Bcanalep.eps} \right\rrangle \oplus \left\llangle \pic{1.2}{0.3}{reid2d.eps} \right\rrangle \right\}, \left\llangle \pic{1.2}{0.3}{reid2c.eps} \right\rrangle \{2\} , 0 , (\partial_{*A}, \partial_{A8} ) \right)
$$
As in the case for oriented links, $C'''$ has null homology and we can conclude showing that the quotient $(C/C')/C'''$ is isomorphic to $\left(\left\llangle \pic{0.9}{0.2}{Bcanalep.eps} \right\rrangle , 0 , 0, 0 \right)$ and using the flatten functor and the lemma. The isomorphism can be found explicitly as in the classic case.

\section{Categorification of the skein relations}

 We note that $\left\llangle \pic{0.5}{0.2}{banp.eps} \right\rrangle = ( 0, W , 0 , 0 ) $, hence $\mathcal{H}^F_*\left(\pic{0.5}{0.2}{banp.eps} \right) = (0, W)$, therefore we have the last equation of the skein relations of the unnormalized Kauffman bracket
$$
\chi_A\left( \mathcal{H}^F_*\left( \pic{0.8}{0.3}{banp.eps} \right) \right) = -A^{-2} - A^2
$$
\begin{teo}
Let $D$ and $D'$ be two link diagrams. Then
$$
\llangle D \sqcup D' \rrangle = \llangle D \rrangle \otimes \llangle D' \rrangle
$$
\begin{proof}
As in the classic case before the addition of the information of the orientation.
\end{proof}
\end{teo}

So we have the second equation of the skein relation, in fact given a diagram $D$
\beq
\chi_A\left( \mathcal{H}^F_* \left( D \sqcup \pic{0.8}{0.3}{banp.eps} \right) \right) & = & \chi_A \left( \left\llangle D \sqcup \pic{0.8}{0.3}{banp.eps} \right\rrangle \right) \\
 & = & \chi_A \left( \llangle D \rrangle \otimes \left\llangle \pic{0.8}{0.3}{banp.eps} \right\rrangle \right) \\
 & = & \chi_A ( \llangle D \rrangle ) \chi_A\left( \left\llangle \pic{0.8}{0.3}{banp.eps} \right\rrangle \right) \\
 & = & \chi_A ( \llangle D \rrangle ) (-A^{-2} - A^2)
\eeq
Let $D$ be a link diagram. We consider one of its crossings $\pic{1.2}{0.3}{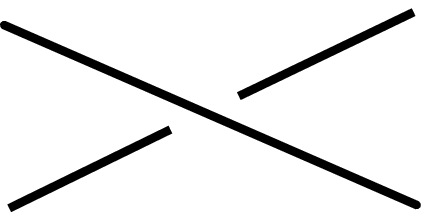}$. We have the following decomposition in direct sum at the level of spaces:
$$
\left\llangle \pic{1.2}{0.3}{incrociop.eps} \right\rrangle = \left\llangle \pic{1.2}{0.3}{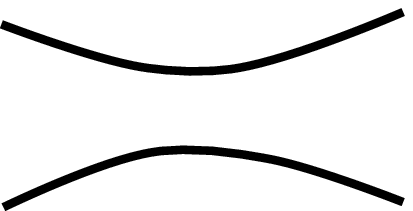} \right\rrangle \{1\} \oplus \left\llangle \pic{1.2}{0.3}{Bcanalep.eps} \right\rrangle \{-1\}
$$
Hence
\beq
\chi_A \left( \mathcal{H}^F_*\left( \pic{1.2}{0.3}{incrociop.eps} \right) \right) & = & \chi_A \left( \left\llangle \pic{1.2}{0.3}{incrociop.eps} \right\rrangle \right) \\
 & = & \chi_A \left( \left\llangle \pic{1.2}{0.3}{Acanalep.eps} \right\rrangle\{1\} \right) + \chi_A \left( \left\llangle \pic{1.2}{0.3}{Bcanalep.eps} \right\rrangle \{-1\} \right) \\
 & = & A\chi_A \left( \mathcal{H}^F_*\left( \pic{1.2}{0.3}{Acanalep.eps} \right) \right) + A^{-1}\chi_A \left( \mathcal{H}^F_*\left( \pic{1.2}{0.3}{Bcanalep.eps} \right) \right)
\eeq
Therefore we have the first equation of the skein relations.
We note that this last equation holds only at the level of spaces and if we want to consider also the differentials we obtain a short exact sequence in $2\!-\!\mathcal{C}_*(\mathcal{G}\textit{r}\mathcal{A}\textit{b})$ where the first map is the inclusion in the first factor and the second is the projection on the second one:
$$
0 \longrightarrow \left\llangle \pic{1.2}{0.3}{Bcanalep.eps} \right\rrangle\{-1\} \longrightarrow \left\llangle \pic{1.2}{0.3}{incrociop.eps} \right\rrangle \longrightarrow \left\llangle \pic{1.2}{0.3}{Acanalep.eps} \right\rrangle\{1\} \longrightarrow 0
$$

So in the end the categorification of the skein relations in homology is:

\begin{enumerate}
\item{
This diagram is exact for each $j\in\Z$
$$
\xymatrix{
\mathcal{H}^F_{0,j+1}\left(\pic{1.2}{0.3}{Bcanalep.eps}\right) \ar[r] & \mathcal{H}_{0,j}\left(\pic{1.2}{0.3}{incrociop.eps}\right) \ar[r] & \mathcal{H}_{0,j-1}\left(\pic{1.2}{0.3}{Acanalep.eps}\right) \ar[d] \\
\mathcal{H}^F_{1,j-1}\left(\pic{1.2}{0.3}{Acanalep.eps}\right) \ar[u] & \mathcal{H}^F_{1,j}\left(\pic{1.2}{0.3}{incrociop.eps}\right) \ar[l] & \mathcal{H}^F_{1,j+1}\left(\pic{1.2}{0.3}{Bcanalep.eps}\right) \ar[l]
}
$$
}

\item{
$$
\mathcal{H}^F_0\left(D \sqcup \pic{1.2}{0.3}{banp.eps}\right) = \mathcal{H}^F_0( D ) \otimes W^\spadesuit
$$
}

\item{
$$
\mathcal{H}^F_*\left(\pic{1.2}{0.3}{banp.eps}\right) = W^\spadesuit
$$
}

\end{enumerate}
The second relation follow from the previous observations and the K\"unneth formula.

\subsection{Positive curls}
Now we see what happens if we follow the proof of the invariance of Khovanov homology for the Reidemeister moves of the first type as presented in \cite{Bar-Natan} and we adapt it.

Given a diagram with a positive curl \pic{0.8}{0.2}{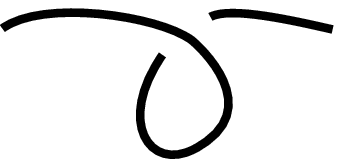} we have the following decomposition in direct sum for the spaces:
$$
\left\llangle \pic{1.0}{0.3}{ricciolopos.eps} \right\rrangle = \left\llangle \pic{1.0}{0.3}{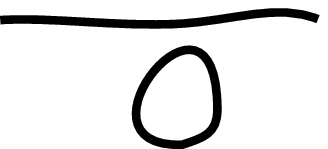} \right\rrangle \{ 1\} \oplus \left\llangle \pic{1.0}{0.3}{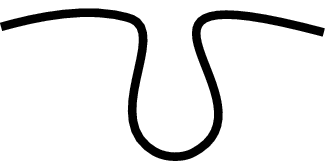} \right\rrangle \{-1\}
$$
By only using the skein relations we have the following result, as in the proposition about the Kauffman bracket to introduce Kauffman's version of the Jones polynomial.
$$
\chi_A \left( \mathcal{H}^F_*\left( \pic{1.0}{0.3}{ricciolopos.eps} \right) \right) = -A^3 \chi_A \left( \mathcal{H}^F_*\left( \pic{1.0}{0.3}{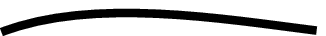} \right) \right)
$$
Following the proof of the invariance of Khovanov homology for the moves of the first type we define a $2$-complex $C$ such that $\textit{Fl}(C)= \left\llangle \pic{0.8}{0.2}{ricciolopos.eps} \right\rrangle$
$$
C := \left( \left\llangle \pic{1.0}{0.3}{riccioloa.eps} \right\rrangle \{1\} , \left\llangle \pic{1.0}{0.3}{ricciolob.eps} \right\rrangle^\spadesuit \{-1\} , \bar m , 0 \right)
$$
Now we define the $2$-subcomplex
$$
C' := \left( \left\llangle \pic{1.0}{0.3}{riccioloa.eps} \right\rrangle_- \{1\} , \left\llangle \pic{1.0}{0.3}{ricciolob.eps} \right\rrangle^\spadesuit \{-1\} , \bar m , 0 \right)
$$
As before $\bar m$ is an isomorphism and hence $C'$ has null homology. $C$ has the same homology as the quotient
$$
\frac{C}{C'} = \left( \left\llangle \pic{1.0}{0.3}{riccioloa.eps} \right\rrangle_+ \{1\}, 0 , 0 , 0 \right)
$$
Therefore using the flatten functor we obtain that
\beq
\mathcal{H}^F_*\left( \pic{1.0}{0.3}{ricciolopos.eps} \right) & \cong & \textit{Fl}( H_*(C/C') ) \\
 & \cong & H_* \left( \left\llangle \pic{1.0}{0.3}{riccioloa.eps} \right\rrangle \{1\} \right) \\
 & \cong & H_* \left( \left( \left\llangle \pic{1.0}{0.3}{riga.eps} \right\rrangle \otimes W_+^\spadesuit \right) \{1\} \right)
\eeq
having identified an graded abelian group $G$ with the $2$-complex $(G,0,0,0)$, and hence $G^\spadesuit = (0, G , 0, 0)$.
\beq
\mathcal{H}^F_*\left( \pic{1.0}{0.3}{ricciolopos.eps} \right) & \cong & H_* \left( \left( \left\llangle \pic{1.0}{0.3}{riga.eps} \right\rrangle \otimes \Z^\spadesuit \{2\} \right) \{1\} \right) \\
 & \cong & H_* \left( \left\llangle \pic{1.0}{0.3}{riga.eps} \right\rrangle^\spadesuit \{3\} \right) \\
 & = & \left( \mathcal{H}^F_* \left( \pic{1.0}{0.3}{riga.eps} \right) \right)^\spadesuit \{3\}
\eeq
Therefore
\beq
\chi_A \left( \mathcal{H}^F_* \left( \pic{1.0}{0.3}{ricciolopos.eps} \right) \right) & = & \chi_A \left( \mathcal{H}^F_* \left( \left( \pic{1.0}{0.3}{riga.eps} \right) \right)^\spadesuit \{3\} \right) \\
 & = & -A^3 \chi_A \left( \mathcal{H}^F_* \left( \pic{1.0}{0.3}{riga.eps} \right) \right)
\eeq
Coherently with what said above.

\section{Returning to the orientations}

As in Kauffman's approach we can consider the oriented link diagrams and add to our construction the information about orientations.

\begin{defn}
Let $A=(A_0,A_1, \partial_0 , \partial_1)$ be a $2$-complex, and $n$ a natural number. We denote with $A^{\spadesuit (n)}$ the $2$-complex obtained reflexing $C$ $n$ times, namely $C^{\spadesuit (n)} = C$ if $n$ is even and $C^{\spadesuit (n)} = C^\spadesuit$ if $n$ is odd.
\end{defn}

\begin{defn}
Given an oriented link diagram $D$ we define the $2$-complex of graded abelian groups
$$
\ddot{C}( D ) := \llangle D \rrangle^{\spadesuit (w(D))} \{-3w(D) \}
$$
where $w(D)$ is the writhe number of $D$.
\end{defn}

\begin{teo}
$H_* \circ \ddot{C} $ up to isomophism is an invariant for oriented links.

\begin{proof}
The operations of shift and reflex commute with the passing in homology. Hence for each link diagram $D$ $H_* ( \ddot{C} (D)) = ( H_* ( \llangle D \rrangle ) )^{\spadesuit(w(D))} \{-3w(D)\} = ( \mathcal{H}^F_*(D))^{\spadesuit(w(D))} \{-3w(D) \}$. Therefore from the invariance of the writhe number and of $\mathcal{H}^F_*$ for the Reidemeister moves of the second and third type and the choice of the order of the crossings, we have the invariance of $H_* \circ \ddot{C}$ for the same modifications. It remain to prove the invariance for the moves of the first type. We consider a diagram with a positive curl \pic{0.8}{0.2}{ricciolopos.eps}. For what said in the previous section we have that
\beq
H_*\left( \ddot{C} \left( \pic{1.0}{0.3}{ricciolopos.eps} \right) \right) & = & \left( \mathcal{H}^F_*\left( \pic{1.0}{0.3}{ricciolopos.eps} \right) \right)^{\spadesuit\left(w \left( \pic{0.5}{0.1}{ricciolopos.eps} \right) \right)} \left\{ -3w\left( \pic{1.0}{0.3}{ricciolopos.eps} \right) \right\} \\
 & \cong & \left( \left( \mathcal{H}^F_* \left( \pic{1.0}{0.3}{riga.eps} \right) \right)^\spadesuit \{3\} \right)^{\spadesuit\left(w \left( \pic{0.5}{0.1}{ricciolopos.eps} \right) \right)} \left\{ -3w\left( \pic{1.0}{0.3}{ricciolopos.eps} \right) \right\} \\
  & \cong & \left( \left( \mathcal{H}^F_* \left( \pic{1.0}{0.3}{riga.eps} \right) \right)^\spadesuit \{3\} \right)^{\spadesuit\left(w \left( \pic{0.5}{0.1}{riga.eps} \right) +1\right)} \left\{ -3w\left( \pic{1.0}{0.3}{riga.eps} \right)  - 3\right\} \\
    & \cong & \left(  \mathcal{H}^F_* \left( \pic{1.0}{0.3}{riga.eps} \right)   \right)^{\spadesuit\left(w \left( \pic{0.5}{0.1}{riga.eps} \right) \right)} \left\{ -3w\left( \pic{1.0}{0.3}{riga.eps} \right)  \right\} \\
 & = & H_*\left( \ddot{C}\left( \pic{1.0}{0.3}{riga.eps} \right) \right)
\eeq
As usual the proof of the invariance for the moves of the first type for the negative curls follows from the invariance for the positive curls and the moves of the second type.
\end{proof}
\end{teo}

\begin{defn}
We denote the composition $H_* \circ \ddot{C}$ with $\ddot{\mathcal{H}}_*$ and in the same way the application defined on the set of the oriented links of $\mathbb{S}^3$.
\end{defn}

\begin{teo}
Let $L$ be an oriented link. Then the graded Euler characteristic of $\ddot{\mathcal{H}}_*(L)$ is equal to the unnormalized Kauffman's version of the Jones polynomial
$$
\chi_A ( \ddot{\mathcal{H}}_*(L) ) = \hat f_L
$$
\begin{proof}
It suffices to prove it at the level of $2$-complexes. Let $D$ be a diagram of $L$.
\beq
\chi_A( \ddot{C}(L) ) & = & \chi_A( (\llangle D \rrangle )^{\spadesuit(w(D) ) } \{ -3w(D)\} ) \\
 & = & (-1)^{w(D)}A^{-3w(D)} \chi_A ( \llangle D \rrangle ) \\
 & = & (-A^3)^{-w(D)} (-A^2-A^{-2}) \langle D \rangle \\
 & = & \hat f_L
\eeq
\end{proof}
\end{teo}

\begin{prop}
Let $D$ be an oriented diagram, $i \in\{0,1\}$ and $j \in 2\Z +1$. Then
$$
\ddot{C}_{i,j}(D) = 0, \ \ddot{\mathcal{H}}_{i,j}(D)= 0
$$
\begin{proof}
The second statement follows from the first one and the arbitrariness of $i\in\{0,1\}$ and $j\in 2\Z +1$. $\ddot{C}_i(D) = ( \bigoplus_{s\ : \ b(s) + |s_A| \in 2\Z + i} W_s(D) )^{\spadesuit(w(D))} \{-3w(D)\} $. Let $s$ be a state of $D$, $W_s(D) \cong   W^{\otimes |s|} \{a(s) -b(s)\}$, where $|s|$ is the number of components of the splitting of $D$ by $s$, $D_s$. By definition of $W$ and of the tensor product of graded abelian groups we have that the only non null homogeneous components of $W^{\otimes |s|}$ are the ones with index that can be obtained summing some copies of $2$ and some copies of $-2$. Hence each odd homogeneous component of $W_s(D)$ is $0$.
\beq
\ddot{C}_{i,j}(D) & = & ( \bigoplus_{s\ : \ b(s) + |s_A| \in 2\Z + i + w(D) } W_s(D) )_{j+3w(D)} \\
 & = & \bigoplus_{s\ : \ b(s) + |s_A| \in 2\Z + i + w(D) } W^{\otimes |s|} )_{j+3w(D) - a(s) + b(s)}
\eeq
$j+3w(D) - a(s) + b(s) = j + 3n_+(D) - 3n_-(D) - n(D) + 2 b(s)$ where $n(D)$ is the number of crossings of $D$ and $n_+(D)$ and $n_-(D)$ are the numbers of positive and negative crossings of $D$. $j+3w(D) - a(s) + b(s) = j + 4 n_+(D) - 2 n_-(D) + 2 b(s)$, since $j$ is odd this number is odd. Hence the component is null.
\end{proof}
\end{prop}

Now we investigate the connections between this new invariant for oriented links, $\ddot{H}_*$, and the classical Khovanov homology, $\mathcal{H}^*$. The first part of the following theorem has been suggested by the famous relation between Khovanov's version of the Jones polynomial and the one of Kauffman: $q= -A^{-2}$.

\begin{teo}\label{teo1}
Let $L$ be an oriented link, $D$ a diagram of $L$ and $N$ the number of components of $L$. Then for each $i\in\{0,1\}$ and $j\in 2\Z$
$$
\begin{array}{cl}
\ddot{H}_{i,j}(L) = \bigoplus_{k \in 2\Z +i} \mathcal{H}^{k,- \frac{j}{2} } (L) & \text{if } j \in 4\Z \\
\ddot{H}_{i,j}(L) = \bigoplus_{k \in 2\Z +i+1} \mathcal{H}^{k,- \frac{j}{2} } (L) & \text{if } j \not\in 4\Z \\
\end{array}
$$
Furthermore for each $i\in\{0,1\}$ and $j\in 2\Z$
$$
\ddot{H}_{i,j} (L) = \bigoplus_{k\in 2\Z +i + n_+(D) + |s_A|} \mathcal{H}^{k, - \frac{j}{2}} (L)
$$
and
\begin{itemize}
\item{if $N\in 2\Z +1$, then $n_+(D) + |s_A| \in 2\Z +1$ and $\ddot{H}_{i,j} (L) = 0$ for each $i\in\{0,1\}$ and $j\in4\Z$;}
\item{if $N\in 2\Z $, then $n_+(D) + |s_A| \in 2\Z $ and $\ddot{H}_{i,j} (L) = 0$ for each $i\in\{0,1\}$ and $j\in 4\Z + 2$.}
\end{itemize}

\begin{proof}
We remind that the classical Khovanov homology is constructed starting from the graded abelian group $V$, that is the free abelian group generated by the elements $v_+$ and $v_-$ where $v_+$ has degree $1$ and $v_-$ has degree $-1$. Let $\alpha : W \rightarrow V$ be the map of abelian groups defined by $\alpha(w_+)= v_-$, $\alpha(w_-) = v_+$. For eache $n\in \mathbb{N}$ we have the map $\alpha^{\otimes n} : W^{\otimes n} \rightarrow V^{\otimes n}$. We note that for eache $j\in 2\Z$ $\alpha((W)_j) = (V)_{-\frac j 2 }$ and $\alpha^{\otimes n}( (W^{\otimes n})_j ) = ( V^{\otimes n} )_{-\frac j 2 }$. From this we can define a map for each state $s$ of $D$ by applying the right shifts: $\alpha_s: W_s(D) \rightarrow V_s(D)$, where $V_s(D) = \left( \bigotimes_{\pic{0.2}{0.1}{banp.eps} \text{ in } D_s} V \right) \{b(s)\}$. For any $j \in 2\Z + n(D)$ $\alpha_s( (W_s(D))_j ) = ( V_s(D) )_{b(s) - \frac{j - a(s) + b(s)}{2} }$. We note that $b(s) - \frac{j - a(s) + b(s)}{2} = - \frac{j-n(D)}{2}$. For each edge $\xi: s \rightarrow s'$ of $D$ the following square in the category of abelian groups (not graded) and morphisms of these, is commutative and the vertical arrows are isomorphisms
$$
\xymatrix{
W_s(D) \ar[d]_{\alpha_s} \ar[r]^{\partial_\xi} & W_{s'}(D) \ar[d]_{\alpha_{s'}} \\
V_s(D) \ar[r]^{d_\xi} & V_{s'}(D)
}
$$
$d_\xi$ is the map of graded abelian groups induced by the edge $\xi$ in the classical Khovanov homology. For each $k\in\Z$ we define the isomorphism of abelian groups $\alpha_k : \bigoplus_{s \ : \ b(s) = k} W_s(D) \rightarrow \llbracket D \rrbracket^k $, where $\llbracket D \rrbracket^k$ is the component in homological degree $k$ of the Khovanov complex of $D$ before the shifts: $\llbracket D \rrbracket^k = \bigoplus_{s \ : \ b(s) = k} V_s(D)$. For any $j\in2\Z+n(K)$, $k\in\Z$ $\alpha_k( (\bigoplus_{s \ : \ b(s) = k} W_s(D) )_j ) = \llbracket D \rrbracket^{k, -\frac{j-n(D)}{2} } $. For any $k\in\Z$ we have the map $\hat d ^k = \sum_{\xi \ :\ |\xi| = k } (-1)^\xi d_\xi :  \llbracket D \rrbracket^k \rightarrow \llbracket D \rrbracket^{k+1}$, and we define $\partial^k := \sum_{\xi \ :\ |\xi| = k} (-1)^\xi \partial_\xi : \bigoplus_{s \ : \ b(s) = k} W_s(D) \rightarrow \bigoplus_{s' \ : \ b(s') = k + 1} W_{s'}(D)$. Hence for each integer $k$ we have the following commutative square in the category of abelian groups with vertical arrows that are isomorphisms
$$
\xymatrix{
\bigoplus_{s \ : \ b(s) = k } W_s(D) \ar[d]_{\alpha_k} \ar[r]^{\partial^k} & \bigoplus_{s' \ : \ b(s') = k + 1} W_{s'}(D) \ar[d]_{\alpha_{k+1}} \\
\llbracket D \rrbracket^k \ar[r]^{\hat d^k} & \llbracket D \rrbracket^{k+1}
}
$$
We note that for each $i\in\{0,1\}$ $\partial_i  = \bigoplus_{k\in 2\Z+i + |s_A|} \partial^k$. Let $i\in\{0,1\}$ and $j\in2\Z$
\beq
\ddot{H}_{i,j}(L) & = & \mathcal{H}^F_{[i+w(D)], j +3w(D)} (D) \\
 & = & \left( \frac{\Ker \partial_{[i+w(D)]}}{ \Im \partial_{[i+w(D) +1]} } \right)_{j+3w(D)} \\
 & \cong & \bigoplus_{k \in 2\Z + i + w(D) + |s_A|} \left( \frac{\Ker \partial^k}{\Im \partial^{k+1}} \right)_{j+3w(D)} \\
 & \cong & \bigoplus_{k \in 2\Z + i + w(D) + |s_A|} \left( \frac{\Ker \hat d^k}{\Im \hat d^{k+1}} \right)_{- \frac{j+3w(D) - n(D)}{2} } \\
 & = & \bigoplus_{k \in 2\Z + i + w(D) + |s_A|} \left( \frac{\Ker \hat d^k}{\Im \hat d^{k+1}} \right)_{- \frac j 2 - n_+(D) + 2n_-(D) } \\
 & = & \bigoplus_{k \in 2\Z + i + n_+(D) + |s_A|} \left( \frac{\Ker d^k}{\Im d^{k+1}} \right)_{- \frac j 2} \\
 & = & \bigoplus_{k \in 2\Z + i + n_+(D) + |s_A|} \mathcal{H}^{k, -\frac j 2} (L)
\eeq
Therefore
\beq
\chi_A (\ddot{H}_*(L)) & = & \sum_{i\in\{0,1\},j\in \Z} (-1)^i A^j \rk \ddot{H}_{i,j} (L) \\
 & = & \sum_{k\in\Z,j\in 2\Z} (-1)^i A^j \rk \mathcal{H}^{k - n_+(D) - |s_A| ,-frac j 2 } (L) \\
 & = & \sum_{k,j\in\Z} (-1)^{n_+(D) + |s_A|} (-1)^i A^{-2j} \rk \mathcal{H}^{i,j}(L) \\
 & = & (-1)^{n_+(D) +|s_A|} \chi_A ( \mathcal{H}^*(L) )(A^{-2})
\eeq
Hence
$$
\hat f_L = (-1)^{n_+(D) + |s_A|} \hat J_L (A^{-2})
$$
where $\hat f$ is the Kauffman's version of the Jones polynomial and $\hat J$ is the Khovanov's one.\\
Let $i\in\{0,1\}$. If $N\in 2\Z +i$, the number of components of $L$, then for each $k \in \Z$, $j\in 2\Z +i +1 $ $\mathcal{H}^{k,j}(L) = 0$. Hence
\begin{itemize}
\item{if $N$ is odd, then $\hat J_L$ has only odd exponents and $\hat f_L$ hasn't exponents multiple of $4$;}
\item{if $N$ is even, then $\hat J_L$ has only exponents even and $\hat f_L$ has only exponents multiple of $4$.}
\end{itemize}
Hence
\beq
\hat f_L & = & \hat J_L (-A^{-2}) \\
 & = & \left\{\begin{array}{cl}
 - \hat J_L (A^{-2}) & \text{if } N \in 2\Z +1 \\
 \hat J_L (A^{-2}) & \text{if } N \in 2\Z
 \end{array}\right. \\
 & = & \left\{\begin{array}{cl}
 - \chi_A( \mathcal{H}^*(L) ) (A^{-2}) & \text{if } N \in 2\Z +1 \\
 \chi_A ( \mathcal{H}^*(L ) ) (A^{-2}) & \text{if } N \in 2\Z
 \end{array}\right. \\
\eeq
It follow that $n_+ (D) + |s_A| $ is congruous modulo $2$ to $N$.\\
Now we suppose that $N$ and $n_+(D) + |s_A|$ are odd. For each $i\in\{0,1\}$ and $j\in 2\Z$ $\ddot{H}_{i,j}(L) = \bigoplus_{k \in2\Z + i +1} \mathcal{H}^{i,-\frac j 2} (L)$. If $j\in 4\Z$ $\ddot{H}_{i,j}(L) = \bigoplus_{k \in2\Z + i +1} \mathcal{H}^{i,-\frac j 2} (L) = 0 = \bigoplus_{k \in2\Z + i } \mathcal{H}^{i,-\frac j 2} (L)$, because $\mathcal{H}^{i,-\frac j 2} (L) = 0$ for each $k$.

If $N$ and $n_+(D) + |s_A|$ are even, for each $i\in\{0,1\}$ and $j\in 2\Z$ $\ddot{H}_{i,j}(L) = \bigoplus_{k \in2\Z + i } \mathcal{H}^{i,-\frac j 2} (L)$. If $j\not\in 4\Z$ $\ddot{H}_{i,j}(L) = \bigoplus_{k \in2\Z + i } \mathcal{H}^{i,-\frac j 2} (L) = 0 = \bigoplus_{k \in2\Z + i +1} \mathcal{H}^{i,-\frac j 2} (L)$, because $\mathcal{H}^{i,-\frac j 2} (L) = 0$ for each $k$.
\end{proof}
\end{teo}

\section{Examples}

$$
\begin{array}{cccccccccccc}
\mathcal{H}^*\left( \pic{1.1}{0.3}{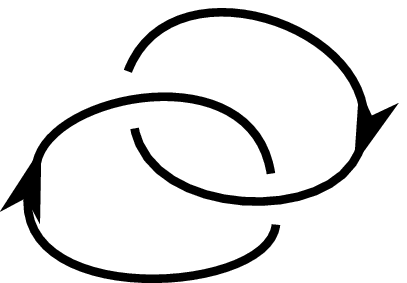} \right) &&&&&&&&&& \\
\mathcal{H}^2 & \ldots & 0 & 0 & 0 & 0 & 0 & \Z & 0 & \Z & 0 &\ldots \\
\mathcal{H}^1 & \ldots & 0 & 0 & 0 & 0 & 0 & 0 & 0 & 0 & 0 & \ldots \\
\mathcal{H}^0 & \ldots & 0 & \Z & 0 & \Z & 0 & 0 & 0 & 0 & 0 & \ldots \\
 & \ldots & -1 & 0 & 1 & 2 & 3 & 4 & 5 & 6 & 7 & \ldots
\end{array}
$$

$$
\begin{array}{cccccccccccc}
\mathcal{H}^F_*\left( \pic{1.1}{0.3}{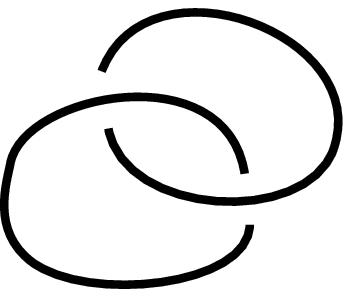} \right) &&&&&&&&&& \\
\mathcal{H}^F_1 & \ldots & 0 & 0 & 0 & 0 & 0 & 0 & \ldots \\
\mathcal{H}^F_0 & \ldots & 0 & \Z & \Z & \Z & \Z & 0 & \ldots \\
 & \ldots & -16 & -12 & -8 & -4 & 0 & 4  & \ldots
\end{array}
$$

$$
\begin{array}{ccccccccccccccc}
\mathcal{H}^*\left( \pic{0.7}{0.2}{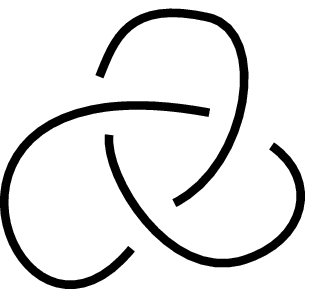} \right)   \\
\mathcal{H}^3 & & \ldots & 0 & 0 & 0 & 0 & 0 & 0 & 0 & (\Z_2)^3 & 0 & \Z & 0 & \ldots \\
\mathcal{H}^2 & & \ldots & 0 & 0 & 0 & 0 & 0 & \Z & 0 & 0 & 0  & 0 & 0 & \ldots \\
\mathcal{H}^1 & & \ldots & 0 & 0 & 0 & 0 & 0 & 0 & 0 & 0 & 0 & 0 & 0 & \ldots \\
\mathcal{H}^0 & & \ldots & 0 & \Z & 0 & \Z & 0 & 0 & 0 & 0 & 0 & 0 & 0 & \ldots \\
\\
 & & \ldots & 0 & 1 & 2 & 3 & 4 & 5 & 6 & 7 & 8 & 9 & 10 & \ldots
\end{array}
$$

$$
\begin{array}{ccccccccccc}
\mathcal{H}^F_*\left( \pic{0.7}{0.2}{trifp.eps} \right)   \\
\mathcal{H}^F_1 & & \ldots & 0 & \Z & (\Z_2)^3 & 0 & 0 & 0 & 0 & \ldots \\
\mathcal{H}^F_0 & & \ldots & 0 & 0 & 0 & \Z & \Z & \Z & 0 & \ldots \\
 & & \ldots & -13 & -9 & -5 & -1 & 3 & 7 & 11 & \ldots
\end{array}
$$

\subsection{Homology is stronger than Kauffman bracket}

Let $\textit{Kh}(L) \in \Z[t,t^{-1},q,q^{-1}]$ be the graded Poincar\'e polynomial of the classical Khovanov homology of the oriented link $L$:
$$
\textit{Kh}(L) := \sum_{i,j\in\Z} t^i q^j \rk \mathcal{H}^{i,j}(L)
$$
Let $\textit{FKh}(L) \in \Z[t,t^{-1},A,A^{-1}]$ be the Poincar\'e polynomial of the homology of the framed link $L$:
$$
\textit{FKh}(L) := \sum_{i \in \{0,1\},j\in\Z} t^i A^j \rk \mathcal{H}^F_{i,j}(L)
$$
Since the theorem \ref{teo1} we know how to obtain $\sum_{i\in\{0,1\}, j\in\Z} t^i A^j \ddot{\mathcal{H}}_{i,j}(L)$ from $\textit{Kh}(L)$ summing some coefficients. If $L$ is represented by the diagram $D$ we also know how to obtain $\mathcal{H}^F_*(D)$ from $\ddot{\mathcal{H}}_*(D)$ using the writhe number, and hence obtain $\textit{FKh}(D)$ from $\textit{Kh}(D)$. In particular if $w(D)=0$ $\mathcal{H}^F_*(D) = \ddot{\mathcal{H}}_*(D)$.

\begin{figure}[htbp]
\begin{center}
\subfigure[$D'_1$]{
$$
\includegraphics[scale=0.5]{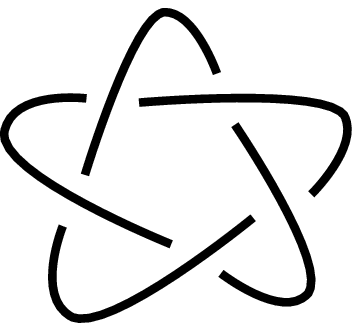}
$$
}\qquad
\subfigure[$D_1$]{
$$
\includegraphics[scale=0.5]{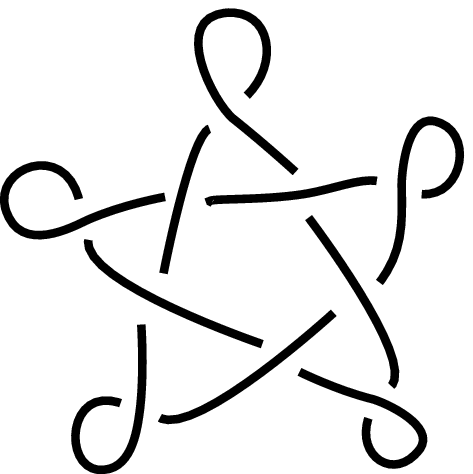}
$$
}
\end{center}
\caption{$5_1$ with two different framings}
\label{figure:5-1}
\end{figure}

\begin{figure}[htbp]
\begin{center}
\subfigure[$D'_2$]{
$$
\includegraphics[scale=0.5]{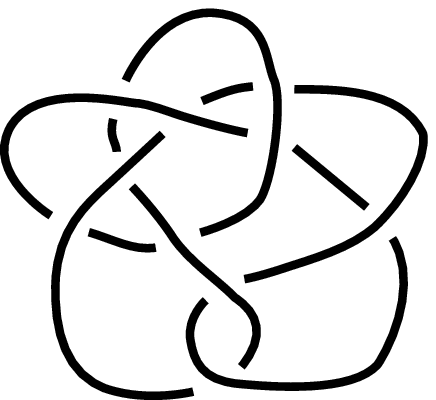}
$$
}\qquad
\subfigure[$D_2$]{
$$
\includegraphics[scale=0.5]{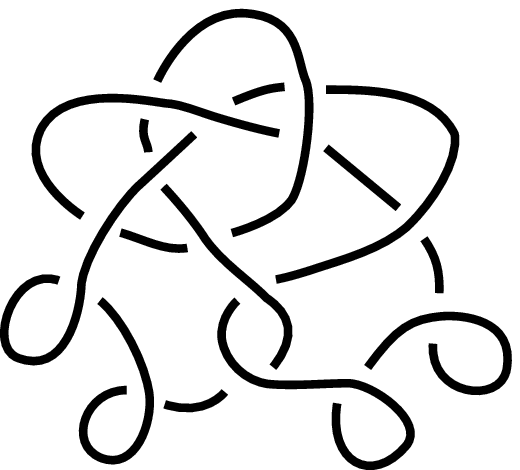}
$$
}
\end{center}
\caption{$10_{132}$ with two different framings}
\label{figure:10-132}
\end{figure}

Let $D'_1$ be the classical diagram of the knot $5_1$ (Figure \ref{figure:5-1} left) and let $D_1$ be the diagram of the knot $5_1$ obtained from $D'_1$ adding $5$ positive curls (Figure \ref{figure:5-1} right). Let $D'_2$ be the classical diagram of the knot $10_{132}$ (Figure \ref{figure:10-132} left) and let $D_2$ be the diagram of $10_{132}$ obtained adding $4$ positive curls (Figure \ref{figure:10-132} right).

Since the writhe number of $D'_1$ is $-5$, the one of $D_1$ is $0$, $w(D'_1)=-5$ and $w(D_1)=0$. And in the same way $w(D'_2)=-4$ and $w(D_2)=0$.

For each diagram $D$ the Kauffman's version of the Jones polynomial of $D$ is $f_D= (-A^3)^{-w(D)} f\langle D \rangle$, therefore
\beq
\langle D_1 \rangle & = & f_{D_1} \\
 & = & f_{D'_1} \\
 & = & -A^{28} + A^{24} - A^{20} + A^{16} + A^8 \\
 & = & f_{D'_2} \\
 & = & f_{D_2} \\
 & = & \langle D_2 \rangle
\eeq
\beq
\textit{Kh}(5_1) & = & q^{-5}+q^{-3}+t^{-5}q^{-15}+t^{-4}q^{-11}+t^{-3}q^{-11}+t^{-2}q^{-7} \\
\textit{Kh}(10_{132}) & = & q^{-3}+q^{-1}+t^{-7}q^{15}+t^{-6}q^{-11}+t^{-5}q^{-11}+t^{-4}q^{-9} \\
 & & +t^{-4}q^{-7}+t^{-3}q^{-9}+t^{-3}q^{-5}+2t^{-2}q^{-5}+t^{-1}q^{-1}
\eeq
\beq
\textit{FKh}(5_1) & = & A^{28} + A^{30} + t( A^6 + A^{10} + A^{14}) + A^{22} \\
\textit{FKh}(10_{132}) & = & A^2 + A^{10} + A^{18} + A^{22} + A^{30}  \\
 & & + t( A^2 + A^6 + 2A^{10} + A^{14} + A^{18} + A^{22} )
\eeq

Hence we have found two framed links with the same Kauffman bracket and different homology.

\section*{Ackoledgements}
The author would like to express his gratitude to Riccardo Benedetti for his useful discussion. He would also like to thank his friends for helping him in the details of this paper.

\end{document}